\documentclass[12pt]{article}
\author{Edward Hanson}
\title{\bf How to recognize a Leonard pair}
\date{}

\usepackage{amssymb}
\usepackage{amsmath}
\usepackage{cite}
\usepackage{verbatim}

\newtheorem{definition}{Definition}[section]
\newtheorem{theorem}[definition]{Theorem}
\newtheorem{proposition}[definition]{Proposition}
\newtheorem{lemma}[definition]{Lemma}

\newtheorem{note}[definition]{Note}

\usepackage[margin=1.0in]{geometry}

\def\fld{\mathbb K}

\begin{document}
\maketitle

\begin{abstract}
\noindent Let $V$ denote a vector space with finite positive dimension. We consider an ordered pair of linear transformations
$A: V\rightarrow V$ and $A^{*}: V\rightarrow V$ that satisfy (i) and (ii) below.
\begin{enumerate}
\item There exists a basis for $V$ with respect to which the matrix representing $A$ is irreducible tridiagonal and the matrix representing $A^{*}$ is diagonal.
\item There exists a basis for $V$ with respect to which the matrix representing $A^{*}$ is irreducible tridiagonal and the matrix representing $A$ is diagonal.
\end{enumerate}
We call such a pair a {\it Leonard pair} on $V$. In the literature, there are some parameters that are used to describe Leonard pairs called the intersection numbers $\{a_{i}\}_{i=0}^{d}$, $\{b_{i}\}_{i=0}^{d-1}$, $\{c_{i}\}_{i=1}^{d}$, and the dual eigenvalues $\{\theta^{*}_{i}\}_{i=0}^{d}$. In this paper, we provide two characterizations of Leonard pairs. For the first characterization, the focus is on the $\{a_{i}\}_{i=0}^{d}$ and $\{\theta^{*}_{i}\}_{i=0}^{d}$. For the second characterization, the focus is on the $\{b_{i}\}_{i=0}^{d-1}$, $\{c_{i}\}_{i=1}^{d}$, and $\{\theta^{*}_{i}\}_{i=0}^{d}$.

\bigskip

\end{abstract}

\section{Introduction} \label{sec:intro}

We begin by recalling the notion of a Leonard pair \cite{T:subconst1, T:Leonard}. We will use the following terms. A square matrix $X$ is called {\it tridiagonal} whenever each nonzero entry lies on either the diagonal, the subdiagonal, or the superdiagonal. Assume $X$ is tridiagonal. Then $X$ is called {\it irreducible} whenever each entry on
the subdiagonal is nonzero and each entry on the superdiagonal is nonzero.

\medskip

\noindent We now define a Leonard pair. For the rest of this paper, $\fld$ will denote a field.

\begin{definition} \label{def:lp} \rm \cite[Definition 1.1]{T:Leonard}
Let $V$ denote a vector space over $\fld$ with finite positive dimension. By a {\it Leonard pair} on $V$, we mean an ordered pair of $\fld$-linear maps $A: V\rightarrow V$ and $A^{*}: V\rightarrow V$ that satisfy (i) and (ii) below.
\begin{enumerate}
\item There exists a basis for $V$ with respect to which the matrix representing $A$ is irreducible tridiagonal and the matrix representing $A^{*}$ is diagonal.
\item There exists a basis for $V$ with respect to which the matrix representing $A^{*}$ is irreducible tridiagonal and the matrix representing $A$ is diagonal.
\end{enumerate}
\end{definition}

\begin{note}
\rm
In a common notational convention, $A^{*}$ denotes the conjugate-transpose of $A$. We are not using this convention. In a Leonard pair $A,A^{*}$, the linear transformations $A$ and $A^{*}$ are arbitrary subject to (i), (ii) above.
\end{note}

\noindent The concept of a Leonard pair originated in the study of $Q$-polynomial distance-regular graphs \cite[p.~260]{BIbook}, \cite[Definition 2.3]{T:subconst1}. Since that time, Leonard pairs have found application in a variety of contexts, such as special functions/orthogonal polynomials \cite{T:Leonard, T:LP24, T:intro, T:TD-D} and representation theory \cite{T:intro, T:madrid}. Motivated by these applications, a number of characterizations of Leonard pairs have been discovered. For instance, there are characterizations of Leonard pairs in terms of orthogonal polynomials \cite[Theorem 19.1]{T:qRacah} \cite[Theorem 4.1]{T:PA}, parameter arrays \cite[Theorem 1.9]{T:Leonard}, upper/lower bidiagonal matrices \cite[Theorem 3.2]{T:PA} \cite[Theorem~17.1]{T:TD-D}, tridiagonal/diagonal matrices \cite[Theorem 25.1]{T:TD-D}, the notion of a tail \cite[Theorem~5.1]{Hanson} \cite[Theorem 10.1]{Hanson3}, and the intersection numbers $\{a_{i}\}_{i=0}^{d}$ \cite[Theorem 5.1]{Hanson2}.

\medskip

\noindent In this paper, we consider the following situation. Fix an integer $d\geq 1$ and consider matrices $A$ and $A^{*}$ over $\fld$ that have the following form:
\begin{equation}
A=\left(
\begin{array}
{ c c c c c c}
  a_{0} & b_{0} &       &       &       & {\bf 0} \\
  c_{1} & a_{1} & b_{1} &       &       & \\
        & c_{2} & \cdot & \cdot &       & \\
        &       & \cdot & \cdot & \cdot & \\
        &       &       & \cdot & \cdot & b_{d-1} \\
{\bf 0} &       &       &       & c_{d} & a_{d}
\end{array}
\right)
\qquad \qquad
A^{*}=\left(
\begin{array}
{ c c c c c c}
\theta^{*}_{0} &                &       &       &       & {\bf 0} \\
               & \theta^{*}_{1} &       &       &       & \\
               &                & \cdot &       &       & \\
               &                &       & \cdot &       & \\
               &                &       &       & \cdot & \\
{\bf 0}        &                &       &       &       & \theta^{*}_{d}
\end{array}
\right). \notag
\end{equation}
It is desirable to have attractive necessary and sufficient conditions for $A, A^{*}$ to form a Leonard pair. In the literature, there exist two kinds of results along this line. For the first kind of result, the focus is on the parameters $\{a_{i}\}_{i=0}^{d}$ and $\{\theta^{*}_{i}\}_{i=0}^{d}$ \cite[Theorem 5.1]{Hanson2}. For the second kind of result, the focus is on the parameters $\{b_{i}\}_{i=0}^{d-1}$, $\{c_{i}\}_{i=1}^{d}$, and $\{\theta^{*}_{i}\}_{i=0}^{d}$ \cite[Theorem~25.1]{T:TD-D}. Each of these results has its drawbacks which we will describe shortly. The present paper has two main theorems, the first of which improves on \cite[Theorem 5.1]{Hanson2} and the second of which improves on \cite[Theorem 25.1]{T:TD-D}. We now describe the drawbacks of \cite[Theorem~5.1]{Hanson2} and \cite[Theorem 25.1]{T:TD-D}, and how our results are an improvement. One shortcoming of \cite[Theorem 5.1]{Hanson2} is that it assumes $A$ is diagonalizable. Our improvement requires no such assumption. The result \cite[Theorem 25.1]{T:TD-D} involves some equations containing the products $\{b_{i-1}c_{i}\}_{i=1}^{d}$, and checking the equations becomes cumbersome. Our improvement avoids this difficulty by treating the $\{b_{i}\}_{i=0}^{d-1}$ and $\{c_{i}\}_{i=1}^{d}$ separately. Our two main results are Theorem~\ref{thm:ai_generalized} and Theorem \ref{thm:main2}.

\section{Leonard systems} \label{sec:LS}

When working with a Leonard pair, it is often convenient to consider a related object called a Leonard system. To prepare for our definition of a Leonard system, we recall a few concepts from linear algebra. From now on, fix an integer $d \geq 0$. Let $\mbox{Mat}_{d+1}(\fld)$ denote the $\fld$-algebra consisting of all $d+1$ by $d+1$ matrices with entries in $\fld$. We index the rows and columns by $0,1,\ldots ,d$. Let $\fld^{d+1}$ denote the vector space over $\fld$ consisting of all $d+1$ by $1$ matrices with entries in $\fld$. We index the rows by $0,1,\ldots ,d$. The algebra $\mbox{Mat}_{d+1}(\fld)$ acts on $\fld^{d+1}$ by left multiplication. Let $V$ denote a vector space over $\fld$ with dimension $d+1$. Let $\mbox{End}(V)$ denote the $\fld$-algebra consisting of the $\fld$-linear maps from $V$ to $V$. The identity of $\mbox{End}(V)$ will be denoted by $I$. The $\fld$-algebra $\mbox{End}(V)$ is isomorphic to $\mbox{Mat}_{d+1}(\fld)$. Let $\{v_{i}\}_{i=0}^{d}$ denote a basis for $V$. For $X\in \mbox{End}(V)$ and $Y\in \mbox{Mat}_{d+1}(\fld)$, we say that $Y$ {\it represents} $X$ {\it with respect to} $\{v_{i}\}_{i=0}^{d}$ whenever $Xv_{j}=\sum_{i=0}^{d}Y_{ij}v_{i}$ for $0\leq j\leq d$. Let $A$ denote an element of $\mbox{End}(V)$. A subspace $W\subseteq V$ will be called an {\it eigenspace} of $A$ whenever $W\neq 0$ and there exists $\theta \in \fld$ such that $W=\{v\in V|Av=\theta v\}$; in this case, $\theta$ is the {\it eigenvalue} of $A$ associated with $W$. We say that $A$ is {\it diagonalizable} whenever $V$ is spanned by the eigenspaces of $A$. We say that $A$ is {\it multiplicity-free} whenever $A$ is diagonalizable and each eigenspace of $A$ has dimension one. By a {\it system of mutually orthogonal idempotents} in $\mbox{End}(V)$, we mean a sequence $\{E_{i}\}_{i=0}^{d}$ of elements in $\mbox{End}(V)$ such that
\begin{equation}
E_{i}E_{j}=\delta_{i,j}E_{i}\qquad \qquad (0 \leq i,j \leq d), \notag
\end{equation}
\begin{equation}
{\rm rank}(E_{i})=1\qquad \qquad (0 \leq i \leq d). \notag
\end{equation}
By a {\it decomposition of $V$}, we mean a sequence $\{U_{i}\}_{i=0}^{d}$ of one-dimensional subspaces of $V$ such that
\begin{equation}
V=\sum_{i=0}^{d}U_{i}\qquad \qquad \text{(direct sum)}. \notag
\end{equation}
The following lemmas are routinely verified.

\begin{lemma} \label{lem:MOidem2}
Let $\{U_{i}\}_{i=0}^{d}$ denote a decomposition of $V$. For $0\leq i\leq d$, define $E_{i}\in \mbox{\rm End}(V)$ such that $(E_{i}-I)U_{i}=0$ and $E_{i}U_{j}=0$ if $j\ne i$ $(0\leq j\leq d)$. Then $\{E_{i}\}_{i=0}^{d}$ is a system of mutually orthogonal idempotents in $\mbox{\rm End}(V)$. Conversely, let $\{E_{i}\}_{i=0}^{d}$ denote a system of mutually orthogonal idempotents in $\mbox{\rm End}(V)$. Define $U_{i}=E_{i}V$ for $0\leq i\leq d$. Then $\{U_{i}\}_{i=0}^{d}$ is a decomposition of $V$.
\end{lemma}

\begin{lemma} \label{lem:EsumI}
Let $\{E_{i}\}_{i=0}^{d}$ denote a system of mutually orthogonal idempotents in $\mbox{\rm End}(V)$. Then $I=\sum_{i=0}^{d}E_{i}$.
\end{lemma}

\noindent Let $A$ denote a multiplicity-free element of $\mbox{End}(V)$ and let $\{\theta_{i}\}^{d}_{i=0}$ denote an ordering of the eigenvalues of $A$. For $0\leq i\leq d$, let $V_{i}$ denote the eigenspace of $A$ for $\theta_{i}$. Then $\{V_{i}\}_{i=0}^{d}$ is a decomposition of $V$; let $\{E_{i}\}_{i=0}^{d}$ denote the corresponding system of mutually orthogonal idempotents from Lemma \ref{lem:MOidem2}. One checks that $A=\sum_{i=0}^{d}\theta_{i}E_{i}$ and $AE_{i}=E_{i}A=\theta_{i}E_{i}$ for $0\leq i\leq d$. Moreover,
\begin{equation}
E_{i}=\prod_{\genfrac{}{}{0pt}{}{0 \leq  j \leq d}{j\not=i}}\frac{A-\theta_{j}I}{\theta_{i}-\theta_{j}}\qquad \qquad (0\leq i\leq d). \notag
\end{equation}
We refer to $E_{i}$ as the {\it primitive idempotent} of $A$ corresponding to $V_{i}$ (or $\theta_{i}$).

\medskip

\noindent We now define a Leonard system.

\begin{definition} \label{def:ls} \rm \cite[Definition 1.4]{T:Leonard}
By a {\it Leonard system} on $V$, we mean a sequence
\begin{equation}
\Phi = (A; \{E_{i}\}_{i=0}^{d}; A^{*}; \{E^{*}_{i}\}_{i=0}^{d}) \notag
\end{equation}
that satisfies (i)--(v) below.
\begin{enumerate}
\item Each of $A,A^{*}$ is a multiplicity-free element of $\mbox{End}(V)$.
\item $\{E_{i}\}_{i=0}^{d}$ is an ordering of the primitive idempotents of $A$.
\item $\{E^{*}_{i}\}_{i=0}^{d}$ is an ordering of the primitive idempotents of $A^{*}$.
\item ${\displaystyle{E^{*}_{i}AE^{*}_{j} =
\begin{cases}
0, & \text{if $\;|i-j|>1$;} \\
\neq 0, & \text{if $\;|i-j|=1$}
\end{cases}
}}
\qquad \qquad (0 \leq i,j\leq d).$
\item ${\displaystyle{E_{i}A^{*}E_{j} =
\begin{cases}
0, & \text{if $\;|i-j|>1$;} \\
\neq 0, & \text{if $\;|i-j|=1$}
\end{cases}
}}
\qquad \qquad (0 \leq i,j\leq d).$
\end{enumerate}
The Leonard system $\Phi$ is said to be {\it over} $\fld$.
\end{definition}

\noindent Let $\Phi = (A; \{E_{i}\}_{i=0}^{d}; A^{*}; \{E^{*}_{i}\}_{i=0}^{d})$ denote a Leonard system on $V$. Then the pair $A, A^{*}$ is a Leonard pair on $V$ said to be {\it associated} with $\Phi$. See \cite[pp. 4--5]{T:Leonard} for the precise connection between Leonard pairs and Leonard systems.

\begin{definition}
\rm
Let $\Phi = (A; \{E_{i}\}_{i=0}^{d}; A^{*}; \{E^{*}_{i}\}_{i=0}^{d})$ denote a Leonard system on $V$. For $0\leq i\leq d$, let $\theta_{i}$ (resp. $\theta^{*}_{i}$) denote the eigenvalue of $A$ (resp. $A^{*}$) associated with $E_{i}V$ (resp. $E^{*}_{i}V$). We call $\{\theta_{i}\}_{i=0}^{d}$ (resp. $\{\theta^{*}_{i}\}_{i=0}^{d}$) the {\it eigenvalue sequence} (resp. {\it dual eigenvalue sequence}) of $\Phi$.
\end{definition}

\begin{definition}
\rm
Let $A, A^{*}$ denote a Leonard pair on $V$. By an {\it eigenvalue sequence} (resp. {\it dual eigenvalue sequence}) of $A, A^{*}$, we mean the eigenvalue sequence (resp. dual eigenvalue sequence) of an associated Leonard system.
\end{definition}

\noindent For the remainder of this section, let $\Phi = (A; \{E_{i}\}_{i=0}^{d}; A^{*}; \{E^{*}_{i}\}_{i=0}^{d})$ denote a Leonard system on $V$ with eigenvalue sequence $\{\theta_{i}\}_{i=0}^{d}$ and dual eigenvalue sequence $\{\theta^{*}_{i}\}_{i=0}^{d}$. To avoid trivialities, we assume $d\geq 1$. By construction, $\{\theta_{i}\}_{i=0}^{d}$ are mutually distinct and contained in $\fld$. Similarly, $\{\theta^{*}_{i}\}_{i=0}^{d}$ are mutually distinct and contained in $\fld$. By \cite[Theorem 12.7]{T:Leonard}, the expressions
\begin{equation} \label{eq:thetarecur}
\frac{\theta_{i-2}-\theta_{i+1}}{\theta_{i-1}-\theta_{i}}, \qquad \qquad \frac{\theta^{*}_{i-2}-\theta^{*}_{i+1}}{\theta^{*}_{i-1}-\theta^{*}_{i}}
\end{equation}
are equal and independent of $i$ for $2\leq i\leq d-1$. Define $\beta\in\fld$ as follows. For $d\geq 3$, let $\beta + 1$ be the common value of (\ref{eq:thetarecur}). For $d\leq 2$, let $\beta$ be arbitrary. By (\ref{eq:thetarecur}), $\theta_{i-1}-\beta\theta_{i}+\theta_{i+1}$ is independent of $i$ for $1\leq i\leq d-1$. Let $\gamma$ denote this common value, so
\begin{equation} \label{eq:gamma}
\theta_{i-1}-\beta\theta_{i}+\theta_{i+1} = \gamma \qquad \qquad (1\leq i\leq d-1).
\end{equation}
For notational convenience, define $\theta_{-1}$ (resp. $\theta_{d+1}$) such that (\ref{eq:gamma}) holds at $i=0$ (resp. $i=d$). Similarly, there exists $\gamma^{*}\in\fld$ such that
\begin{equation} \label{eq:gamma*}
\theta^{*}_{i-1}-\beta\theta^{*}_{i}+\theta^{*}_{i+1} = \gamma^{*} \qquad \qquad (1\leq i\leq d-1).
\end{equation}
For notational convenience, define $\theta^{*}_{-1}$ (resp. $\theta^{*}_{d+1}$) such that (\ref{eq:gamma*}) holds at $i=0$ (resp. $i=d$). Choose $0\neq u\in E_{0}V$. By \cite[Lemma 5.1]{T:LP24}, $E^{*}_{i}u$ is a basis for $E^{*}_{i}V$ $(0\leq i\leq d)$. Moreover, $\{E^{*}_{i}u\}_{i=0}^{d}$ is a basis for $V$. By Lemma \ref{lem:EsumI},
\begin{equation} \label{eq:usum}
u = \sum_{i=0}^{d} E^{*}_{i}u.
\end{equation}
With respect to the basis $\{E^{*}_{i}u\}_{i=0}^{d}$, the matrices representing $A$ and $A^{*}$ take the form
\begin{equation}
A:\left(
\begin{array}
{ c c c c c c}
  a_{0} & b_{0} &       &       &       & {\bf 0} \\
  c_{1} & a_{1} & b_{1} &       &       & \\
        & c_{2} & \cdot & \cdot &       & \\
        &       & \cdot & \cdot & \cdot & \\
        &       &       & \cdot & \cdot & b_{d-1} \\
{\bf 0} &       &       &       & c_{d} & a_{d}
\end{array}
\right)
\qquad \qquad
A^{*}:\left(
\begin{array}
{ c c c c c c}
\theta^{*}_{0} &                &       &       &       & {\bf 0} \\
               & \theta^{*}_{1} &       &       &       & \\
               &                & \cdot &       &       & \\
               &                &       & \cdot &       & \\
               &                &       &       & \cdot & \\
{\bf 0}        &                &       &       &       & \theta^{*}_{d}
\end{array}
\right), \notag
\end{equation}
for some scalars $a_{i}, b_{i}, c_{i} \in \fld$ with $c_{i}b_{i-1}\neq 0$ for $1\leq i\leq d$. We call the scalars $\{a_{i}\}_{i=0}^{d}$, $\{b_{i}\}_{i=0}^{d-1}$, and $\{c_{i}\}_{i=1}^{d}$ the \emph{intersection numbers} of $\Phi$. By (\ref{eq:usum}) and since $Au = \theta_{0}u$,
\begin{equation} \label{eq:CRS}
c_{i} + a_{i} + b_{i} = \theta_{0} \qquad \qquad (0\leq i\leq d),
\end{equation}
where $c_{0} = b_{d} = 0$. By \cite[Definition 7.1 and Lemma 7.2]{T:qRacah},
\begin{equation}
a_{i}=\mbox{tr}(E^{*}_{i}A) \qquad \qquad (0 \leq i \leq d), \notag
\end{equation}
where $\mbox{tr}$ denotes trace. The next equation involves the intersection number $a_{0}^{*}$ for the Leonard system $(A^{*}; \{E^{*}_{i}\}_{i=0}^{d}; A; \{E_{i}\}_{i=0}^{d})$. By \cite[Lemma 9.2]{T:PA},
\begin{equation} \label{eq:bici}
c_{i}(\theta^{*}_{i-1}-\theta^{*}_{i})-b_{i}(\theta^{*}_{i}-\theta^{*}_{i+1}) = (\theta_{1}-\theta_{0})(\theta^{*}_{i}-a^{*}_{0}) \qquad \qquad  (0\leq i\leq d).
\end{equation}
By \cite[Theorem 5.3]{TV}, there exist $\omega, \eta^{*} \in \fld$ such that
\begin{equation} \label{eq:ai_cond1}
a_{i}(\theta^{*}_{i}-\theta^{*}_{i-1})(\theta^{*}_{i}-\theta^{*}_{i+1})=\gamma \theta^{*2}_{i}+\omega \theta^{*}_{i}+\eta^{*} \qquad \qquad (0\leq i\leq d).
\end{equation}
Using (\ref{eq:ai_cond1}), we obtain
\begin{equation} \label{eq:ai_cond2}
a_{i}(\theta^{*}_{i}-\theta^{*}_{i+1})+a_{i-1}(\theta^{*}_{i-1}-\theta^{*}_{i-2})-\gamma(\theta^{*}_{i-1}+\theta^{*}_{i})=\omega \qquad \qquad (1\leq i\leq d).
\end{equation}
\begin{proposition}
With the above notation,
\begin{equation} \label{eq:Omega}
c_{i}(\theta^{*}_{i-1}-\theta^{*}_{i+1})-b_{i-1}(\theta^{*}_{i-2}-\theta^{*}_{i})-(\theta_{0}-\theta_{-1})(\theta^{*}_{i-1}+\theta^{*}_{i}) = \Omega \qquad \qquad (1\leq i\leq d),
\end{equation}
where $\Omega = 2\theta_{0}(a^{*}_{0}-\gamma^{*})-2\theta_{1}a^{*}_{0}-\omega$.
\end{proposition}
\noindent {\it Proof:} Let the integer $i$ be given. In (\ref{eq:bici}), eliminate $c_{i}$ using (\ref{eq:CRS}) to obtain
\begin{equation} \label{eq:bi_expr}
b_{i}(\theta^{*}_{i+1}-\theta^{*}_{i-1}) = \theta_{1}(\theta^{*}_{i}-a^{*}_{0})+\theta_{0}(a^{*}_{0}-\theta^{*}_{i-1})-a_{i}(\theta^{*}_{i}-\theta^{*}_{i-1}).
\end{equation}
In (\ref{eq:bi_expr}), replace $i$ by $i-1$ to obtain
\begin{equation} \label{eq:bi-1_expr}
b_{i-1}(\theta^{*}_{i}-\theta^{*}_{i-2}) = \theta_{1}(\theta^{*}_{i-1}-a^{*}_{0})+\theta_{0}(a^{*}_{0}-\theta^{*}_{i-2})-a_{i-1}(\theta^{*}_{i-1}-\theta^{*}_{i-2}).
\end{equation}
In (\ref{eq:bici}), eliminate $b_{i}$ using (\ref{eq:CRS}) to obtain
\begin{equation} \label{eq:ci_expr}
c_{i}(\theta^{*}_{i-1}-\theta^{*}_{i+1}) = \theta_{1}(\theta^{*}_{i}-a^{*}_{0})+\theta_{0}(a^{*}_{0}-\theta^{*}_{i+1})-a_{i}(\theta^{*}_{i}-\theta^{*}_{i+1}).
\end{equation}
Adding (\ref{eq:bi-1_expr}) to (\ref{eq:ci_expr}) and simplifying the result using (\ref{eq:gamma}), (\ref{eq:gamma*}), and (\ref{eq:ai_cond2}), we routinely obtain (\ref{eq:Omega}). \hfill $\Box$ \\

\section{The first main theorem}

In this section, we obtain our first main result. It involves the following setup. Fix an integer $d\geq 1$. Let $V$ denote a vector space over $\fld$ of dimension $d+1$. Let $\{E^{*}_{i}\}^{d}_{i=0}$ denote a system of mutually orthogonal idempotents in $\mbox{End}(V)$. Define $A\in\mbox{End}(V)$ such that
\begin{equation} \label{eq:E*AE*}
{ \displaystyle E^{*}_{i}AE^{*}_{j} =
\begin{cases}
0, & \text{if $\;|i-j|>1$;} \\
\neq 0, & \text{if $\;|i-j|=1$}
\end{cases}
}
\qquad \qquad (0 \leq i,j\leq d).
\end{equation}
Let $\{\theta^{*}_{i}\}^{d}_{i=0}$ denote scalars in $\fld$ and define
\begin{equation} \label{eq:A*sum}
A^{*}=\sum_{i=0}^{d}\theta^{*}_{i}E^{*}_{i}.
\end{equation}
Let $\{\theta_{i}\}^{d}_{i=0}$ denote any scalars in $\fld$.

\begin{theorem} \label{thm:ai_generalized}
With the above notation, suppose the following {\rm (i)--(vi)} hold.
\begin{enumerate}
\item[\rm (i)] $\theta_{i}\neq \theta_{j}$ if $i\neq j$ $(0\leq i,j\leq d)$.
\item[\rm (ii)] $\theta^{*}_{i}\neq \theta^{*}_{0}$ $(1\leq i \leq d)$.
\item[\rm (iii)] There exist $\beta, \gamma^{*} \in \fld$ such that
\begin{equation} \label{eq:TTR0}
\gamma^{*}=\theta^{*}_{i-1}-\beta \theta^{*}_{i}+\theta^{*}_{i+1} \qquad \qquad (1\leq i \leq d-1).
\end{equation}
Define $\theta^{*}_{-1}$ {\rm (}resp. $\theta^{*}_{d+1}${\rm )} such that {\rm (\ref{eq:TTR0})} holds at $i=0$ {\rm (}resp. $i=d${\rm )}.
\item[\rm (iv)] There exist nonzero vectors $v_{0}, v_{1}\in V$ such that
\begin{equation} \label{eq:v0v1}
Av_{0} = \theta_{0}v_{0}, \qquad \qquad Av_{1} = \theta_{1}v_{1}, \qquad \qquad A^{*}v_{0}-v_{1}\in\fld v_{0}.
\end{equation}
\item[\rm (v)] There exist $\gamma, \omega, \eta^{*} \in \fld$ such that for $0\leq i\leq d$,
\begin{equation} \label{eq:ai_assum}
a_{i}(\theta^{*}_{i}-\theta^{*}_{i-1})(\theta^{*}_{i}-\theta^{*}_{i+1})=\gamma \theta^{*2}_{i}+\omega \theta^{*}_{i}+\eta^{*},
\end{equation}
where $a_{i}=\mbox{\rm tr}(E^{*}_{i}A)$.
\item[\rm (vi)] $\theta_{i-1}-\beta\theta_{i}+\theta_{i+1}=\gamma$ $(1\leq i\leq d-1)$.
\end{enumerate}
Then $A, A^{*}$ is a Leonard pair on $V$ with eigenvalue sequence $\{\theta_{i}\}_{i=0}^{d}$ and dual eigenvalue sequence $\{\theta^{*}_{i}\}_{i=0}^{d}$.
\end{theorem}
\noindent {\it Proof:} By (\ref{eq:E*AE*}) and \cite[Corollary 3.4]{Hanson}, the elements $A$ and $E^{*}_{0}$ together generate $\mbox{End}(V)$. Using (\ref{eq:A*sum}) and the fact that $\{ E^{*}_{i}\}_{i=0}^{d}$ are mutually orthogonal idempotents, we obtain
\begin{equation}
E^{*}_{0} = \prod_{j=1}^{d}\frac{A^{*}-\theta^{*}_{j}I}{\theta^{*}_{0}-\theta^{*}_{j}}. \notag
\end{equation}
Consequently, $\mbox{End}(V)$ is generated by $A$ and $A^{*}$. The vector space $V$ is irreducible as an $\mbox{End}(V)$-module, so $V$ is irreducible as a module for $A, A^{*}$.

\medskip

\noindent By \cite[Lemma 3.5]{Hanson}, there exists a unique antiautomorphism $\dagger$ of $\mbox{End}(V)$ such that $A^{\dagger} = A$ and $E^{*\dagger}_{i} = E^{*}_{i}$ for $0\leq i\leq d$. By this and (\ref{eq:A*sum}), $A^{*\dagger} = A^{*}$.

\medskip

\noindent Recall the scalars $\theta^{*}_{-1}$ and $\theta^{*}_{d+1}$ from below (\ref{eq:TTR0}). By construction,
\begin{equation} \label{eq:TTR2}
\gamma^{*}=\theta^{*}_{i-1}-\beta \theta^{*}_{i}+\theta^{*}_{i+1} \qquad \qquad (0\leq i \leq d).
\end{equation}
We claim that the scalar
\begin{equation} \label{eq:delta*}
\theta ^{*2}_{i-1}-\beta \theta ^{*}_{i-1}\theta ^{*}_{i}+\theta ^{*2}_{i}-\gamma ^{*}(\theta ^{*}_{i-1}+\theta ^{*}_{i})
\end{equation}
is independent of $i$ for $0\leq i\leq d+1$. Denote this scalar by $p_{i}$. For $0\leq i\leq d$,
\begin{equation}
p_{i}-p_{i+1}=(\theta^{*}_{i-1}-\theta^{*}_{i+1})(\theta^{*}_{i-1}-\beta\theta^{*}_{i}+\theta^{*}_{i+1}-\gamma^{*}). \notag
\end{equation}
In this equation, the right-hand side equals $0$ by (\ref{eq:TTR2}). Consequently, $p_{i}$ is independent of $i$ for $0\leq i\leq d+1$. The claim is now proven. Let $\delta^{*}$ denote the common value of (\ref{eq:delta*}), so
\begin{equation} \label{eq:delta*2}
\theta ^{*2}_{i-1}-\beta \theta ^{*}_{i-1}\theta ^{*}_{i}+\theta ^{*2}_{i}-\gamma ^{*}(\theta ^{*}_{i-1}+\theta ^{*}_{i}) = \delta^{*} \qquad \qquad (0\leq i\leq d+1).
\end{equation}

\medskip

\noindent We now show that
\begin{equation} \label{eq:theta*prod}
(\theta^{*}_{i}-\theta^{*}_{i-1})(\theta^{*}_{i}-\theta^{*}_{i+1})=(2-\beta)\theta^{*2}_{i}-2\gamma^{*}\theta^{*}_{i}-\delta^{*} \qquad \qquad (0\leq i \leq d).
\end{equation}
To verify (\ref{eq:theta*prod}), in the right-hand side, replace $\delta^{*}$ by (\ref{eq:delta*}) and eliminate both occurrences of $\gamma^{*}$ in the resulting expression using (\ref{eq:TTR2}). We have now verified (\ref{eq:theta*prod}).

\medskip

\noindent For notational convenience, we introduce a $2$-variable polynomial
\begin{equation} \label{eq:P}
P(\lambda,\mu)=\lambda^{2}-\beta \lambda \mu+\mu^{2}-\gamma^{*}(\lambda +\mu)-\delta^{*}.
\end{equation}
We now claim that
\begin{equation} \label{eq:AW2}
A^{*2}A-\beta A^{*}AA^{*}+AA^{*2}-\gamma^{*}(AA^{*}+A^{*}A)-\delta^{*}A=\gamma A^{*2}+\omega A^{*}+\eta^{*}I.
\end{equation}
In (\ref{eq:AW2}), let $C$ denote the left-hand side minus the right-hand side. We show $C=0$. Using $I=\sum_{i=0}^{d}E^{*}_{i}$, we obtain
\begin{align}
C &= (E^{*}_{0}+E^{*}_{1}+\cdots +E^{*}_{d})C(E^{*}_{0}+E^{*}_{1}+\cdots +E^{*}_{d}) \notag \\
  &= \sum_{i=0}^{d}\sum_{j=0}^{d} E^{*}_{i}CE^{*}_{j}. \notag
\end{align}
For $0\leq i,j\leq d$, we show $E^{*}_{i}CE^{*}_{j}=0$. Using $E^{*}_{i}A^{*}=\theta^{*}_{i}E^{*}_{i}$ and $A^{*}E^{*}_{j}=\theta^{*}_{j}E^{*}_{j}$,
\begin{equation} \label{eq:E*iCE*j}
E^{*}_{i}CE^{*}_{j}=E^{*}_{i}AE^{*}_{j} P(\theta^{*}_{i},\theta^{*}_{j})-\delta_{i,j}(\gamma \theta^{*2}_{i}+\omega \theta^{*}_{i}+\eta^{*})E^{*}_{i}.
\end{equation}
To further examine (\ref{eq:E*iCE*j}), we consider two cases. First assume $i\neq j$. In this case, (\ref{eq:E*iCE*j}) becomes
\begin{equation}
E^{*}_{i}CE^{*}_{j}=E^{*}_{i}AE^{*}_{j} P(\theta^{*}_{i},\theta^{*}_{j}). \notag
\end{equation}
If $|i-j|>1$, then $E^{*}_{i}AE^{*}_{j}=0$ by (\ref{eq:E*AE*}). If $|i-j|=1$, then $P(\theta^{*}_{i},\theta^{*}_{j})=0$ by (\ref{eq:delta*2}). Therefore, $E^{*}_{i}CE^{*}_{j}=0$ under our present assumption that $i\neq j$. Next assume $i=j$. In this case, (\ref{eq:E*iCE*j}) becomes
\begin{equation} \label{eq:E*iCE*i}
E^{*}_{i}CE^{*}_{i}=E^{*}_{i}AE^{*}_{i}P(\theta^{*}_{i},\theta^{*}_{i})-(\gamma \theta^{*2}_{i}+\omega \theta^{*}_{i}+\eta^{*})E^{*}_{i}.
\end{equation}
By (\ref{eq:theta*prod}) and (\ref{eq:P}), we find $P(\theta^{*}_{i},\theta^{*}_{i})=(\theta^{*}_{i}-\theta^{*}_{i-1})(\theta^{*}_{i}-\theta^{*}_{i+1})$. By \cite[Proposition 3.6]{Hanson2}, $E^{*}_{i}AE^{*}_{i}=a_{i}E^{*}_{i}$. Evaluating the right-hand side of (\ref{eq:E*iCE*i}) using these comments, we find that it equals $E^{*}_{i}$ times
\begin{equation} \label{eq:aithetas}
a_{i}(\theta^{*}_{i}-\theta^{*}_{i-1})(\theta^{*}_{i}-\theta^{*}_{i+1})-\gamma \theta^{*2}_{i}-\omega \theta^{*}_{i}-\eta^{*}.
\end{equation}
The scalar (\ref{eq:aithetas}) is equal to 0 by (\ref{eq:ai_assum}), so $E^{*}_{i}CE^{*}_{i}=0$. We have now shown $E^{*}_{i}CE^{*}_{j}=0$ for $0\leq i,j\leq d$. Therefore, $C=0$. We have now verified (\ref{eq:AW2}).

\medskip

\noindent We now claim that for $1\leq i\leq d$, there exists a nonzero vector $v_{i}\in V$ such that both
\begin{equation} \label{eq:v_i_claim}
Av_{i}=\theta_{i}v_{i}, \qquad \qquad A^{*}v_{i-1}-v_{i}\in \mbox{span}(v_{0},\ldots,v_{i-1}),
\end{equation}
where $v_{0}$ is from (\ref{eq:v0v1}). We prove the claim by induction on $i$. The case $i=1$ follows by condition (iv). Next assume $i\geq 2$. Note that $v_{0}, v_{1}, \ldots, v_{i-1}$ are linearly independent, because they are eigenvectors for $A$ with distinct eigenvalues. For $0\leq j\leq i-1$, define $W_{j} = \mbox{span}(v_{0},\ldots ,v_{j})$. By construction,
\begin{equation} \label{eq:Wchain}
W_{0}\subseteq W_{1}\subseteq\cdots\subseteq W_{i-1}.
\end{equation}
By induction,
\begin{align}
    AW_{j} &\subseteq W_{j} \qquad \qquad (0\leq j\leq i-1), \label{eq:AW} \\
A^{*}W_{j} &\subseteq W_{j+1} \qquad \qquad (0\leq j\leq i-2). \label{eq:A*W}
\end{align}
We apply both sides of (\ref{eq:AW2}) to $v_{i-2}$ and evaluate the result using $Av_{i-2} = \theta_{i-2}v_{i-2}$. This gives
\begin{equation} \label{eq:AW2_cons}
(A+\theta_{i-2}-\gamma)A^{*2}v_{i-2} - (\gamma^{*}A+\theta_{i-2}\gamma^{*}+\omega)A^{*}v_{i-2} - \beta A^{*}AA^{*}v_{i-2} - (\delta^{*}\theta_{i-2}+\eta^{*})v_{i-2} = 0.
\end{equation}
For notational convenience, define
\begin{equation} \label{eq:w_i-2}
w_{i-2} = A^{*}v_{i-2}-v_{i-1}.
\end{equation}
Evaluate (\ref{eq:AW2_cons}) using (\ref{eq:w_i-2}), and simplify the result using $Av_{i-1} = \theta_{i-1}v_{i-1}$ and $\beta\theta_{i-1} = \theta_{i-2}+\theta_{i}-\gamma$. This gives
\begin{equation}
\begin{split}
&(A-\theta_{i})A^{*}v_{i-1} + (A+\theta_{i-2}-\gamma )A^{*}w_{i-2} - (\gamma^{*}\theta_{i-1}+\gamma^{*}\theta_{i-2}+\omega)v_{i-1} \\
&- \beta A^{*}Aw_{i-2} - (\gamma^{*}A+\gamma^{*}\theta_{i-2}+\omega)w_{i-2} - (\delta^{*}\theta_{i-2}+\eta^{*})v_{i-2} = 0. \label{eq:AW2_cons2}
\end{split}
\end{equation}
By (\ref{eq:v_i_claim}), (\ref{eq:w_i-2}), and induction, $w_{i-2}\in W_{i-2}$. Using (\ref{eq:AW}) and (\ref{eq:A*W}), $Aw_{i-2}\in W_{i-2}$ and $A^{*}Aw_{i-2}\in W_{i-1}$. Using these comments to simplify (\ref{eq:AW2_cons2}), we obtain
\begin{equation} \label{eq:A*v_i-1}
(A-\theta_{i})A^{*}v_{i-1}\in W_{i-1}.
\end{equation}
We now show that $A^{*}v_{i-1}\not\in W_{i-1}$. Suppose $A^{*}v_{i-1}\in W_{i-1}$. By this, together with (\ref{eq:Wchain}) and (\ref{eq:A*W}), $A^{*}W_{i-1}\subseteq W_{i-1}$. By (\ref{eq:AW}), $AW_{i-1}\subseteq W_{i-1}$. Comparing the dimensions of $W_{i-1}$ and $V$, we obtain $W_{i-1}\neq V$. This contradicts the fact that $V$ is irreducible as a module for $A, A^{*}$. We have shown that $A^{*}v_{i-1}\not\in W_{i-1}$. Let $H$ denote the subspace of $V$ spanned by $W_{i-1}$ and $A^{*}v_{i-1}$. The vectors $v_{0}, \ldots, v_{i-1}, A^{*}v_{i-1}$ form a basis for $H$. Recall that $Av_{j} = \theta_{j}v_{j}$ for $0\leq j\leq i-1$. By this and (\ref{eq:A*v_i-1}), $AH\subseteq H$ and the action of $A$ on $H$ has characteristic polynomial $(\lambda-\theta_{0})(\lambda-\theta_{1})\cdots(\lambda-\theta_{i})$. By condition (i), the roots of this characteristic polynomial are mutually distinct, so $A$ is diagonalizable on $H$ with eigenvalues $\theta_{0}, \ldots, \theta_{i}$. Let $0\neq v_{i}\in H$ denote an eigenvector for $A$ with eigenvalue $\theta_{i}$. So $Av_{i} = \theta_{i}v_{i}$. Note that $v_{i}\not\in W_{i-1}$, so there exists $0\neq\epsilon\in\fld$ such that $A^{*}v_{i-1}-\epsilon v_{i}\in W_{i-1}$. Replacing $v_{i}$ by $\epsilon v_{i}$, we may assume $\epsilon = 1$. We have shown $A^{*}v_{i-1}-v_{i}\in\mbox{span}(v_{0},\ldots ,v_{i-1})$. The claim is proven.

\medskip

\noindent By construction and since $\{\theta_{i}\}_{i=0}^{d}$ are mutually distinct, $\{v_{i}\}_{i=0}^{d}$ is a basis for $V$ consisting of eigenvectors for $A$. It follows that $A$ is multiplicity-free. For $0\leq i\leq d$, let $E_{i}$ denote the primitive idempotent of $A$ corresponding to $\theta_{i}$. We now show that $\Phi = (A; \{E_{i}\}_{i=0}^{d}; A^{*}; \{E^{*}_{i}\}_{i=0}^{d})$ is a Leonard system on $V$. To do this, we verify conditions (i)--(v) of Definition \ref{def:ls}. Definition \ref{def:ls}(ii) holds by construction and Definition \ref{def:ls}(iv) holds by (\ref{eq:E*AE*}). It is convenient to check the remaining conditions in a nonstandard order. Consider Definition \ref{def:ls}(v). By (\ref{eq:v_i_claim}),
\begin{equation} \label{eq:EA*E1}
{ \displaystyle E_{i}A^{*}E_{j} =
\begin{cases}
0, & \text{if $\;i-j>1$;} \\
\neq 0, & \text{if $\;i-j=1$}
\end{cases}
}
\qquad \qquad (0 \leq i,j\leq d).
\end{equation}
Applying $\dagger$,
\begin{equation} \label{eq:EA*E2}
{ \displaystyle E_{i}A^{*}E_{j} =
\begin{cases}
0, & \text{if $\;j-i>1$;} \\
\neq 0, & \text{if $\;j-i=1$}
\end{cases}
}
\qquad \qquad (0 \leq i,j\leq d).
\end{equation}
Definition \ref{def:ls}(v) holds by (\ref{eq:EA*E1}) and (\ref{eq:EA*E2}). To obtain Definition \ref{def:ls}(i), we show that $A^{*}$ is multiplicity-free. The map $A^{*}$ is given in (\ref{eq:A*sum}). By assumption, $\{E^{*}_{i}\}_{i=0}^{d}$ are mutually orthogonal idempotents in $\mbox{End}(V)$. Therefore, by Lemma \ref{lem:MOidem2}, the sum $V = \sum_{i=0}^{d}E^{*}_{i}V$ is direct and $E^{*}_{i}V$ has dimension $1$ for $0\leq i\leq d$. By (\ref{eq:A*sum}), $(A^{*}-\theta^{*}_{i})E_{i}V = 0$ for $0\leq i\leq d$. By these comments, $A^{*}$ is diagonalizable. To show that $A^{*}$ is multiplicity-free, we show that $\{\theta^{*}_{i}\}_{i=0}^{d}$ are mutually distinct. Define a polynomial $\psi(\lambda)=\prod_{i=0}^{d}(\lambda -\theta^{*}_{i})$ and note that $\psi(A^{*})=0$. The elements $\{A^{*i}\}_{i=0}^{d}$ are linearly independent by Definition \ref{def:ls}(v) and \cite[Lemma 3.1]{Hanson}, so the minimal polynomial of $A^{*}$ has degree $d+1$. Therefore, the minimal polynomial of $A^{*}$ is precisely $\psi(\lambda)$. Because $A^{*}$ is diagonalizable, the roots $\{\theta^{*}_{i}\}_{i=0}^{d}$ of $\psi(\lambda)$ are mutually distinct. Therefore, $A^{*}$ is multiplicity-free as desired. We have established Definition \ref{def:ls}(i). By (\ref{eq:A*sum}) and since $A^{*}$ is multiplicity-free, we see that $\{E^{*}_{i}\}_{i=0}^{d}$ is an ordering of the primitive idempotents of $A^{*}$. This gives Definition \ref{def:ls}(iii). By these comments, $\Phi$ is a Leonard system on $V$. Consequently, $A, A^{*}$ is a Leonard pair on $V$ with eigenvalue sequence $\{\theta_{i}\}_{i=0}^{d}$ and dual eigenvalue sequence $\{\theta^{*}_{i}\}_{i=0}^{d}$. \hfill $\Box$ \\

\section{The second main theorem}

In this section, we obtain our second main result.

\begin{theorem} \label{thm:main2}
Fix an integer $d\geq 1$. Suppose there exist scalars $\{\theta_{i}\}_{i=0}^{d}$, $\{\theta_{i}^{*}\}_{i=0}^{d}$, and $\{a_{i}\}_{i=0}^{d}$, $\{b_{i}\}_{i=0}^{d-1}$, $\{c_{i}\}_{i=1}^{d}$ in $\fld$ such that the following {\rm (i)--(viii)} hold.
\begin{enumerate}
\item[\rm (i)] $\theta_{i}\neq \theta_{j}$ if $i\neq j$ $(0\leq i, j\leq d)$.
\item[\rm (ii)] $\theta^{*}_{i}\neq \theta^{*}_{0}$ $(1\leq i\leq d)$.
\item[\rm (iii)] There exist $\beta, \gamma^{*} \in \fld$ such that
\begin{equation} \label{eq:TTR1}
\gamma^{*}=\theta^{*}_{i-1}-\beta\theta^{*}_{i}+\theta^{*}_{i+1} \qquad \qquad (1\leq i \leq d-1).
\end{equation}
Define $\theta^{*}_{-1}$ {\rm (}resp. $\theta^{*}_{d+1}${\rm )} such that {\rm (\ref{eq:TTR1})} holds at $i=0$ {\rm (}resp. $i=d${\rm )}.
\item[\rm (iv)] $b_{i-1}c_{i}\neq 0$ for $1\leq i\leq d$.
\item[\rm (v)] $c_{i}+a_{i}+b_{i}=\theta_{0}$ for $0\leq i \leq d$, where $b_{d}=c_{0}=0$.
\item[\rm (vi)] There exists $a^{*}_{0}\in \fld$ such that
\begin{equation} \label{eq:bici2}
c_{i}(\theta^{*}_{i-1}-\theta^{*}_{i})-b_{i}(\theta^{*}_{i}-\theta^{*}_{i+1}) = (\theta_{1}-\theta_{0})(\theta^{*}_{i}-a^{*}_{0}) \qquad \qquad (0\leq i \leq d).
\end{equation}
\item[\rm (vii)] There exists $\theta_{-1}\in \fld$ such that
\begin{equation} \label{eq:bi-1_ci2}
c_{i}(\theta^{*}_{i-1}-\theta^{*}_{i+1})-b_{i-1}(\theta^{*}_{i-2}-\theta^{*}_{i})-(\theta_{0}-\theta_{-1})(\theta^{*}_{i-1}+\theta^{*}_{i})
\end{equation}
is independent of $i$ for $1\leq i\leq d$.
\item[\rm (viii)] Define $\gamma=\theta_{-1}-\beta\theta_{0}+\theta_{1}$. Then
\begin{equation} \label{eq:gamma2}
\theta_{i-1}-\beta\theta_{i}+\theta_{i+1}=\gamma \qquad \qquad (1\leq i\leq d-1).
\end{equation}
\end{enumerate}
Then there exists a Leonard system over $\fld$ with eigenvalue sequence $\{\theta_{i}\}_{i=0}^{d}$, dual eigenvalue sequence $\{\theta^{*}_{i}\}_{i=0}^{d}$, and intersection numbers $\{a_{i}\}_{i=0}^{d}$, $\{b_{i}\}_{i=0}^{d-1}$, $\{c_{i}\}_{i=1}^{d}$.
\end{theorem}
\noindent {\it Proof:} Define the vector space $V = \fld^{d+1}$. We identify $\mbox{End}(V)$ with $\mbox{Mat}_{d+1}(\fld)$. Define $A, A^{*}\in\mbox{Mat}_{d+1}(\fld)$ as follows:
\begin{equation} \label{eq:A_A*}
A=\left(
\begin{array}
{ c c c c c c}
  a_{0} & b_{0} &       &       &       & {\bf 0} \\
  c_{1} & a_{1} & b_{1} &       &       & \\
        & c_{2} & \cdot & \cdot &       & \\
        &       & \cdot & \cdot & \cdot & \\
        &       &       & \cdot & \cdot & b_{d-1} \\
{\bf 0} &       &       &       & c_{d} & a_{d}
\end{array}
\right)
\qquad \qquad
A^{*}=\left(
\begin{array}
{ c c c c c c}
\theta^{*}_{0} &                &       &       &       & {\bf 0} \\
               & \theta^{*}_{1} &       &       &       & \\
               &                & \cdot &       &       & \\
               &                &       & \cdot &       & \\
               &                &       &       & \cdot & \\
{\bf 0}        &                &       &       &       & \theta^{*}_{d}
\end{array}
\right).
\end{equation}
For $0\leq i\leq d$, define $E^{*}_{i}\in\mbox{Mat}_{d+1}(\fld)$ with $(i,i)$-entry $1$ and all other entries $0$. The elements $\{E^{*}_{i}\}_{i=0}^{d}$ are mutually orthogonal idempotents. Note that $A, A^{*}$ and $\{E^{*}_{i}\}_{i=0}^{d}$ satisfy the conditions stated above Theorem \ref{thm:ai_generalized}.

\medskip

\noindent We now show that $A, A^{*}$ is a Leonard pair. Our strategy is to invoke Theorem \ref{thm:ai_generalized}. We now check the conditions of Theorem \ref{thm:ai_generalized}. First note that Theorem \ref{thm:ai_generalized}(i), Theorem \ref{thm:ai_generalized}(ii), and Theorem \ref{thm:ai_generalized}(iii) are satisfied by conditions (i). (ii), and (iii) in the present theorem, respectively. We now verify Theorem \ref{thm:ai_generalized}(iv). Let $v_{0}\in V$ denote the vector with every component equal to $1$. By condition (v) in the present theorem, $Av_{0} = \theta_{0}v_{0}$. Combining conditions (v) and (vi) in the present theorem, we obtain
\begin{equation} \label{eq:stailTTR}
c_{i}(\theta^{*}_{i-1}-a^{*}_{0})+a_{i}(\theta^{*}_{i}-a^{*}_{0})+b_{i}(\theta^{*}_{i+1}-a^{*}_{0}) = \theta_{1}(\theta^{*}_{i}-a^{*}_{0}) \qquad \qquad (0\leq i \leq d).
\end{equation}
Let $v_{1}\in V$ denote the vector with $i^{\text{th}}$ component $\theta^{*}_{i}-a^{*}_{0}$ for $0\leq i\leq d$. By condition (ii) in the present theorem, $v_{1} \neq 0$. By (\ref{eq:A_A*}) and (\ref{eq:stailTTR}), we obtain $Av_{1} = \theta_{1}v_{1}$ and $A^{*}v_{0}-v_{1} = a^{*}_{0}v_{0}$. This implies Theorem \ref{thm:ai_generalized}(iv).

\medskip

\noindent We now show Theorem \ref{thm:ai_generalized}(v). Evaluating (\ref{eq:bici2}) using condition (v), we obtain
\begin{equation} \label{eq:aibi_expr}
(\theta_{0}-a_{i})(\theta^{*}_{i-1}-\theta^{*}_{i})-b_{i}(\theta^{*}_{i-1}-\theta^{*}_{i+1}) = (\theta_{1}-\theta_{0})(\theta^{*}_{i}-a^{*}_{0}) \qquad \qquad  (0\leq i\leq d).
\end{equation}
Rearranging the terms in (\ref{eq:aibi_expr}), we obtain
\begin{equation} \label{eq:bi_expr2}
b_{i}(\theta^{*}_{i+1}-\theta^{*}_{i-1}) = \theta_{1}(\theta^{*}_{i}-a^{*}_{0})+\theta_{0}(a^{*}_{0}-\theta^{*}_{i-1})-a_{i}(\theta^{*}_{i}-\theta^{*}_{i-1}) \qquad \qquad  (0\leq i\leq d).
\end{equation}
Evaluating (\ref{eq:bici2}) using condition (v), we similarly obtain
\begin{equation} \label{eq:ci_expr2}
c_{i}(\theta^{*}_{i-1}-\theta^{*}_{i+1}) = \theta_{1}(\theta^{*}_{i}-a^{*}_{0})+\theta_{0}(a^{*}_{0}-\theta^{*}_{i+1})-a_{i}(\theta^{*}_{i}-\theta^{*}_{i+1}) \qquad \qquad  (0\leq i\leq d).
\end{equation}
For $1\leq i\leq d$, consider the equation obtained from (\ref{eq:bi_expr2}) by replacing $i$ with $i-1$. Add this to (\ref{eq:ci_expr2}) to obtain
\begin{equation}
\begin{split}
&c_{i}(\theta^{*}_{i-1}-\theta^{*}_{i+1}) + b_{i-1}(\theta^{*}_{i}-\theta^{*}_{i-2}) \\
&= \theta_{1}(\theta^{*}_{i}+\theta^{*}_{i-1}-2a^{*}_{0})+\theta_{0}(2a^{*}_{0}-\theta^{*}_{i+1}-\theta^{*}_{i-2})-a_{i}(\theta^{*}_{i}-\theta^{*}_{i+1})-a_{i-1}(\theta^{*}_{i-1}-\theta^{*}_{i-2}) \label{eq:bi-1_ci_pre2}
\end{split}
\end{equation}
for $1\leq i\leq d$.

\medskip

\noindent Let $\Omega$ denote the common value of (\ref{eq:bi-1_ci2}). By (\ref{eq:bi-1_ci_pre2}) and condition (vii) in the present theorem,
\begin{equation}
\begin{split}
\theta_{1}(\theta^{*}_{i}+\theta^{*}_{i-1}-2a^{*}_{0})+\theta_{0}(2a^{*}_{0}-\theta^{*}_{i+1}-\theta^{*}_{i}-\theta^{*}_{i-1}-\theta^{*}_{i-2}) \\
-a_{i}(\theta^{*}_{i}-\theta^{*}_{i+1})-a_{i-1}(\theta^{*}_{i-1}-\theta^{*}_{i-2})+\theta_{-1}(\theta^{*}_{i-1}+\theta^{*}_{i}) = \Omega \label{eq:Omega2}
\end{split}
\end{equation}
for $1\leq i\leq d$. In (\ref{eq:Omega2}), eliminate $\theta_{-1}$ using $\gamma = \theta_{-1}-\beta\theta_{0}+\theta_{1}$. Evaluating the results using (\ref{eq:TTR1}), we obtain
\begin{equation} \label{eq:omega}
a_{i}(\theta^{*}_{i}-\theta^{*}_{i+1}) + a_{i-1}(\theta^{*}_{i-1}-\theta^{*}_{i-2}) - \gamma(\theta^{*}_{i-1}+\theta^{*}_{i}) = 2\theta_{0}(a^{*}_{0}-\gamma^{*}) - 2\theta_{1}a^{*}_{0} - \Omega
\end{equation}
for $1\leq i\leq d$. Let $\omega$ denote the right-hand side of (\ref{eq:omega}). So,
\begin{equation} \label{eq:ai_cond2b}
a_{i}(\theta^{*}_{i}-\theta^{*}_{i+1})+a_{i-1}(\theta^{*}_{i-1}-\theta^{*}_{i-2})-\gamma(\theta^{*}_{i-1}+\theta^{*}_{i})=\omega \qquad \qquad (1\leq i\leq d).
\end{equation}
For $1\leq i\leq d$, we multiply each side of (\ref{eq:ai_cond2b}) by $\theta^{*}_{i}-\theta^{*}_{i-1}$. After some rearranging, we obtain
\begin{equation}
a_{i}(\theta^{*}_{i}-\theta^{*}_{i-1})(\theta^{*}_{i}-\theta^{*}_{i+1})-\gamma\theta^{*2}_{i}-\omega\theta^{*}_{i} = a_{i-1}(\theta^{*}_{i-1}-\theta^{*}_{i-2})(\theta^{*}_{i-1}-\theta^{*}_{i})-\gamma\theta^{*2}_{i-1}-\omega\theta^{*}_{i-1} \qquad (1\leq i\leq d). \notag
\end{equation}
Consequently, the scalar
\begin{equation} \label{eq:eta*}
a_{i}(\theta^{*}_{i}-\theta^{*}_{i-1})(\theta^{*}_{i}-\theta^{*}_{i+1})-\gamma\theta^{*2}_{i}-\omega\theta^{*}_{i}
\end{equation}
is independent of $i$ for $0\leq i\leq d$. Let $\eta^{*}$ denote the common value of (\ref{eq:eta*}). So,
\begin{equation}
a_{i}(\theta^{*}_{i}-\theta^{*}_{i-1})(\theta^{*}_{i}-\theta^{*}_{i+1})=\gamma \theta^{*2}_{i}+\omega \theta^{*}_{i}+\eta^{*} \qquad \qquad (0\leq i\leq d). \notag
\end{equation}
By the equation on the left in (\ref{eq:A_A*}) and by the definition of $E^{*}_{i}$ following (\ref{eq:A_A*}), we routinely obtain $a_{i} = \mbox{tr}(E^{*}_{i}A)$ $(0\leq i\leq d)$. This establishes Theorem \ref{thm:ai_generalized}(v). Theorem \ref{thm:ai_generalized}(vi) follows from (\ref{eq:gamma2}). We have established the conditions of Theorem \ref{thm:ai_generalized}. Therefore, the pair $A, A^{*}$ is a Leonard pair on $V$ with eigenvalue sequence $\{\theta_{i}\}_{i=0}^{d}$ and dual eigenvalue sequence $\{\theta^{*}_{i}\}_{i=0}^{d}$. For $0\leq i\leq d$, $E^{*}_{i}$ is the primitive idempotent of $A^{*}$ associated with $\theta^{*}_{i}$. For $0\leq i\leq d$, let $E_{i}$ denote the primitive idempotent of $A$ associated with the eigenvalue $\theta_{i}$. By construction, the sequence $\Phi = (A; \{E_{i}\}^{d}_{i=0}; A^{*}; \{E^{*}_{i}\}^{d}_{i=0})$ is a Leonard system on $V$ with eigenvalue sequence $\{\theta_{i}\}_{i=0}^{d}$ and dual eigenvalue sequence $\{\theta^{*}_{i}\}_{i=0}^{d}$. By the equation on the left in (\ref{eq:A_A*}), $\Phi$ has intersection numbers $\{a_{i}\}_{i=0}^{d}$, $\{b_{i}\}_{i=0}^{d-1}$, and $\{c_{i}\}_{i=1}^{d}$. \hfill $\Box$ \\

\section{Three applications of Theorem \ref{thm:main2}}

In this section, we illustrate Theorem \ref{thm:main2} with three examples.

\begin{proposition} \label{prop:Krawtchouk}
Fix an integer $d\geq 1$. Assume that the characteristic of $\fld$ is zero or an odd prime greater than $d$. Define
\begin{align}
    \theta_{i} &= d-2i \qquad \qquad (0\leq i\leq d), \label{eq:Ktheta} \\
\theta^{*}_{i} &= d-2i \qquad \qquad (0\leq i\leq d), \label{eq:Ktheta*} \\
         b_{i} &= d-i \qquad \qquad (0\leq i\leq d-1), \label{eq:Kbi} \\
         c_{i} &= i \qquad \qquad (1\leq i\leq d), \label{eq:Kci} \\
         a_{i} &= 0 \qquad \qquad (0\leq i\leq d).
\end{align}
Then the conditions of Theorem \ref{thm:main2} are satisfied with
\begin{align}
 \beta = 2, \qquad \qquad \gamma = 0, \qquad \qquad \gamma^{*} = 0, \qquad \qquad \theta_{-1} &= d+2, \\
 \theta^{*}_{-1} = d+2, \qquad \qquad \theta^{*}_{d+1} = -d-2, \qquad \qquad a^{*}_{0} &= 0. \label{eq:Kextra2}
\end{align}
\end{proposition}
\noindent {\it Proof:} Using the data (\ref{eq:Ktheta})--(\ref{eq:Kextra2}), one routinely verifies that each of conditions {\rm (i)--(viii)} from Theorem \ref{thm:main2} holds. \hfill $\Box$ \\

\begin{note}
\rm
Referring to Proposition \ref{prop:Krawtchouk}, the corresponding Leonard system from Theorem~\ref{thm:main2} is said to have Krawtchouk type; see \cite[Section 24]{T:madrid}.
\end{note}

\begin{proposition} \label{prop:q-Racah}
Let $\fld$ be arbitrary and fix an integer $d\geq 1$. Let $a$, $b$, $c$, and $q$ denote nonzero scalars in $\fld$ such that each of the following hold.
\begin{itemize}
\item $q^{2i}\neq 1$ for $1\leq i\leq d$.
\item Neither of $a^{2}, b^{2}$ is among $q^{2d-2}, q^{2d-4}, \ldots , q^{2-2d}$.
\item None of $abc, a^{-1}bc, ab^{-1}c, abc^{-1}$ is among $q^{d-1}, q^{d-3}, \ldots , q^{1-d}$.
\end{itemize}
Define
\begin{align}
    \theta_{i} &= aq^{2i-d}+a^{-1}q^{d-2i} \qquad \qquad (0\leq i\leq d), \label{eq:qrtheta} \\
\theta^{*}_{i} &= bq^{2i-d}+b^{-1}q^{d-2i} \qquad \qquad (0\leq i\leq d), \label{eq:qrtheta*} \\
         b_{0} &= \frac{(q^{d}-q^{-d})(cq-a^{-1}b^{-1}q^{d})(q^{-1}-abc^{-1}q^{-d})}{bq^{1-d}-b^{-1}q^{d-1}}, \\
         b_{i} &= \frac{(q^{d-i}-q^{i-d})(bq^{i-d}-b^{-1}q^{d-i})(cq^{i+1}-a^{-1}b^{-1}q^{d-i})(q^{-i-1}-abc^{-1}q^{i-d})}{(bq^{2i-d}-b^{-1}q^{d-2i})(bq^{2i-d+1}-b^{-1}q^{d-2i-1})} \qquad (1\leq i\leq d-1), \label{eq:qrbi1} \\
         c_{i} &= \frac{(q^{i}-q^{-i})(bq^{i}-b^{-1}q^{-i})(a^{-1}q^{i-d-1}-b^{-1}c^{-1}q^{-i})(bq^{i}-acq^{d-i+1})}{(bq^{2i-d-1}-b^{-1}q^{d-2i+1})(bq^{2i-d}-b^{-1}q^{d-2i})} \qquad \qquad (1\leq i\leq d-1), \label{eq:qrci1} \\
         c_{d} &= \frac{(q^{d}-q^{-d})(a^{-1}q^{-1}-b^{-1}c^{-1}q^{-d})(bq^{d}-acq)}{bq^{d-1}-b^{-1}q^{1-d}}, \label{eq:qrcd1} \\
         a_{i} &= \theta_{0}-b_{i}-c_{i} \qquad \qquad (0\leq i\leq d),
\end{align}
where $b_{d} = c_{0} = 0$. Then the conditions of Theorem \ref{thm:main2} are satisfied with
\begin{align}
 \beta &= q^{2}+q^{-2}, \qquad \qquad \gamma = 0, \qquad \qquad \gamma^{*} = 0, \qquad \qquad \theta_{-1} = aq^{-d-2}+a^{-1}q^{d+2}, \\
 \theta^{*}_{-1} &= bq^{-d-2}+b^{-1}q^{d+2}, \qquad \qquad \theta^{*}_{d+1} = bq^{d+2}+b^{-1}q^{-d-2}, \\
 a^{*}_{0} &= \frac{(b+b^{-1})(aq-a^{-1}q^{-1})-(c+c^{-1})(q^{d}-q^{-d})}{aq^{1-d}-a^{-1}q^{d-1}}. \label{eq:qra0*}
\end{align}
\end{proposition}
\noindent {\it Proof:} Using the data (\ref{eq:qrtheta})--(\ref{eq:qra0*}), one routinely verifies that each of conditions {\rm (i)--(viii)} from Theorem \ref{thm:main2} holds. In this calculation, it is useful to note that
\begin{align}
        \theta_{i}-\theta_{j} &= (aq^{i+j-d}-a^{-1}q^{d-i-j})(q^{i-j}-q^{j-i}), \notag \\
\theta^{*}_{i}-\theta^{*}_{j} &= (bq^{i+j-d}-b^{-1}q^{d-i-j})(q^{i-j}-q^{j-i}), \notag
\end{align}
for $0\leq i,j\leq d$. In Theorem \ref{thm:main2}(vii), expression (\ref{eq:bi-1_ci2}) is equal to
\begin{equation}
(q^{2}-q^{-2})\left((q^{d+1}-q^{-d-1})(c+c^{-1})-(a-a^{-1})(b+b^{-1})\right) \notag
\end{equation}
for $1\leq i\leq d$. \hfill $\Box$ \\

\begin{note}
\rm
Referring to Proposition \ref{prop:q-Racah}, the corresponding Leonard system from Theorem~\ref{thm:main2} is said to have $q$-Racah type; see \cite[Section 5]{Huang}.
\end{note}

\noindent In applications, we are often presented with a tridiagonal matrix and a diagonal matrix, each with numerical entries, and we wish to know whether this is a Leonard pair. In our next example, we illustrate how to proceed using Theorem \ref{thm:main2}.

\begin{proposition} \label{prop:KT}
Assume that the characteristic of $\fld$ is zero. Define $d=5$ and
\begin{alignat}{11}
\theta^{*}_{0} &=& 3, \qquad \theta^{*}_{1} &=& \frac{93}{35}, \qquad \theta^{*}_{2} &=& \frac{69}{35}, \qquad \theta^{*}_{3} &=& \frac{33}{35}, \qquad \theta^{*}_{4} &=& -\frac{3}{7},  \qquad \theta^{*}_{5} &=& -\frac{15}{7}, \label{eq:KTtheta*} \\
         b_{0} &=& 3, \qquad b_{1} &=& \frac{64}{35}, \qquad b_{2} &=& \frac{243}{175}, \qquad b_{3} &=& \frac{48}{49}, \qquad b_{4} &=& \frac{11}{21}, \qquad \quad & & \label{eq:KTbi} \\
         \quad & & \qquad c_{1} &=& 1, \qquad c_{2} &=& \frac{192}{175}, \qquad c_{3} &=& \frac{243}{245}, \qquad c_{4} &=& \frac{16}{21}, \qquad c_{5} &=& \frac{3}{7}, \label{eq:KTci} \\
         a_{0} &=& 0, \qquad a_{1} &=& \frac{6}{35}, \qquad a_{2} &=& \frac{18}{35}, \qquad a_{3} &=& \frac{36}{35}, \qquad a_{4} &=& \frac{12}{7}, \qquad a_{5} &=& \frac{18}{7}. \label{eq:KTai}
\end{alignat}
Then the conditions of Theorem \ref{thm:main2} are satisfied with
\begin{align}
 \theta_{i} &= \theta^{*}_{i} \qquad (0\leq i\leq 5), \notag \\
 \beta &= 2, \qquad \gamma = \gamma^{*} = -\frac{12}{35}, \qquad \theta_{-1} = \theta^{*}_{-1} = 3, \qquad \theta^{*}_{6} = -\frac{21}{5}, \qquad a^{*}_{0} = 0. \notag
\end{align}
\end{proposition}
\noindent {\it Proof:} We now verify conditions {\rm (i)--(viii)} in Theorem \ref{thm:main2}. Theorem \ref{thm:main2}(ii) holds by (\ref{eq:KTtheta*}). Theorem \ref{thm:main2}(iv) holds by (\ref{eq:KTbi}) and (\ref{eq:KTci}). Concerning Theorem \ref{thm:main2}(iii), using the data (\ref{eq:KTtheta*}), we evaluate (\ref{eq:TTR1}) at $i = 1$ and $i = 2$ to compute $\beta$ and $\gamma^{*}$. We then verify (\ref{eq:TTR1}) and compute $\theta^{*}_{-1}$ and $\theta^{*}_{6}$. We have now verified Theorem \ref{thm:main2}(iii). Using the data (\ref{eq:KTbi})--(\ref{eq:KTai}), we verify that Theorem \ref{thm:main2}(v) holds with $\theta_{0} = 3 = \theta^{*}_{0}$. Concerning Theorem \ref{thm:main2}(vi), using the data (\ref{eq:KTtheta*})--(\ref{eq:KTci}), we evaluate (\ref{eq:bici2}) at $i = 0$ and $i = 1$. We routinely solve for $\theta_{1}$ and $a^{*}_{0}$, and verify (\ref{eq:bici2}). We have now verified Theorem \ref{thm:main2}(vi). Concerning Theorem \ref{thm:main2}(vii), using the data (\ref{eq:KTtheta*})--(\ref{eq:KTci}), we evaluate (\ref{eq:bi-1_ci2}) at $i=1$ to obtain $\theta_{-1}$ and, using that value, we routinely verify (\ref{eq:bi-1_ci2}). We have now verified Theorem \ref{thm:main2}(vii). We obtain $\gamma$ using the first equation in Theorem \ref{thm:main2}(viii). We define $\theta_{2}$, $\theta_{3}$, $\theta_{4}$, and $\theta_{5}$ so that (\ref{eq:gamma2}) holds. We obtain $\theta_{i} = \theta^{*}_{i}$ ($0\leq i\leq 5$). Note that Theorem \ref{thm:main2}(i) is satisfied. We have now verified each of conditions {\rm (i)--(viii)} from Theorem \ref{thm:main2}. \hfill $\Box$ \\

\begin{note}
\rm
Referring to Proposition \ref{prop:KT}, the corresponding Leonard system from Theorem~\ref{thm:main2} is said to have Racah type; see \cite[Example 35.9]{T:madrid}.
\end{note}

\section{The first and second split sequence}

Consider the Leonard system from Definition \ref{def:ls}. In \cite{T:Leonard}, this Leonard system was described using a sequence of scalars called its parameter array. A parameter array takes the form $\left(\{\theta_{i}\}_{i=0}^{d}, \{\theta^{*}_{i}\}_{i=0}^{d}, \{\varphi_{i}\}_{i=1}^{d}, \{\phi_{i}\}_{i=1}^{d}\right)$, where $\{\theta_{i}\}_{i=0}^{d}$ is the eigenvalue sequence and $\{\theta^{*}_{i}\}_{i=0}^{d}$ is the dual eigenvalue sequence. The sequences $\{\varphi_{i}\}_{i=1}^{d}$ and $\{\phi_{i}\}_{i=1}^{d}$ are called the {\it first and second split sequences}, respectively \cite[p. 5]{T:PA}. It follows from \cite[Definition 23.1 and Theorem 23.5]{T:madrid} that for $d\geq 1$,
\begin{align}
\varphi_{1} &= b_{0}(\theta^{*}_{1}-\theta^{*}_{0}), \label{eq:varphi1} \\
\varphi_{i} &= b_{i-1}\frac{(\theta^{*}_{i}-\theta^{*}_{0})\cdots (\theta^{*}_{i}-\theta^{*}_{i-1})}{(\theta^{*}_{i-1}-\theta^{*}_{0})\cdots (\theta^{*}_{i-1}-\theta^{*}_{i-2})} \qquad \qquad (2\leq i\leq d), \label{eq:varphii} \\
\phi_{i} &= c_{i}\frac{(\theta^{*}_{i-1}-\theta^{*}_{d})\cdots (\theta^{*}_{i-1}-\theta^{*}_{i})}{(\theta^{*}_{i}-\theta^{*}_{d})\cdots (\theta^{*}_{i}-\theta^{*}_{i+1})} \qquad \qquad (1\leq i\leq d-1), \label{eq:phii} \\
\phi_{d} &= c_{d}(\theta^{*}_{d-1}-\theta^{*}_{d}). \label{eq:phid}
\end{align}

\medskip

\noindent Assume that our Leonard system is the one from Proposition \ref{prop:Krawtchouk}. Using (\ref{eq:Ktheta*})--(\ref{eq:Kci}) to simplify (\ref{eq:varphi1})--(\ref{eq:phid}), we obtain
\begin{align}
\varphi_{i} &= -2i(d-i+1) \qquad \qquad (1\leq i\leq d), \notag \\
   \phi_{i} &= 2i(d-i+1) \qquad \qquad (1\leq i\leq d). \notag
\end{align}
This matches the data presented in \cite[Section 16]{T:LP24}.

\medskip

\noindent Next, assume that our Leonard system is the one from Proposition \ref{prop:q-Racah}. Using (\ref{eq:qrtheta*})--(\ref{eq:qrcd1}) to simplify (\ref{eq:varphi1})--(\ref{eq:phid}), we find that for $1\leq i\leq d$,
\begin{align}
\varphi_{i} &= a^{-1}b^{-1}q^{d+1}(q^{i}-q^{-i})(q^{i-d-1}-q^{d-i+1})(q^{-i}-abcq^{i-d-1})(q^{-i}-abc^{-1}q^{i-d-1}), \notag \\
\phi_{i} &= ab^{-1}q^{d+1}(q^{i}-q^{-i})(q^{i-d-1}-q^{d-i+1})(q^{-i}-a^{-1}bcq^{i-d-1})(q^{-i}-a^{-1}bc^{-1}q^{i-d-1}). \notag
\end{align}
This matches the data presented in \cite[Definition 6.1]{Huang}.

\medskip

\noindent Finally, assume that our Leonard system is the one from Proposition \ref{prop:KT}. Using (\ref{eq:KTtheta*})--(\ref{eq:KTci}) to simplify (\ref{eq:varphi1})--(\ref{eq:phid}), we find that
\begin{alignat}{11}
\varphi_{1} &=& -\frac{36}{35}, \qquad \varphi_{2} &=& -\frac{4608}{1225}, \qquad \varphi_{3} &=& -\frac{8748}{1225}, \qquad \varphi_{4} &=& -\frac{2304}{245}, \qquad \varphi_{5} &=& -\frac{396}{49}, \notag \\
\phi_{1} &=& \frac{36}{49}, \qquad \phi_{2} &=& \frac{2304}{1225}, \qquad \phi_{3} &=& \frac{2916}{1225}, \qquad \phi_{4} &=& \frac{2304}{1225}, \qquad \phi_{5} &=& \frac{36}{49}. \notag
\end{alignat}

\section{Acknowledgment}

The author thanks Paul Terwilliger for providing many valuable ideas and detailed suggestions, and John Caughman for providing the data for the example from Proposition \ref{prop:KT}.

\bigskip

\noindent Edward Hanson \hfil\break
\noindent Department of Mathematics \hfil\break
\noindent University of Wisconsin \hfil\break
\noindent 480 Lincoln Drive \hfil\break
\noindent Madison, WI 53706-1388 USA \hfil\break
\noindent email: {\tt hanson@math.wisc.edu }\hfil\break

\end{document}